\documentclass[12pt]{article}
\usepackage{amsthm,amssymb,amscd,amsmath,latexsym,indentfirst,color,amsfonts}
\usepackage[mathcal]{eucal}
\theoremstyle{plain}
\newtheorem{theorem}{Theorem}[section]
\newtheorem{prop}[theorem]{Proposition}
\newtheorem{coro}[theorem]{Corollary}
\newtheorem{lemma}[theorem]{Lemma}
\newtheorem{ex}[theorem]{Example}

\newtheorem{rem}[theorem]{Remark}

\newtheorem{defi}[theorem]{Definition}

\newcommand{\ble}{\begin {lemma}}
\newcommand{\ele}{\end {lemma}}
\newcommand{\bprp}{\begin {proposition}}
\newcommand{\eprp}{\end {proposition}}
\newcommand{\bthm}{\begin {theorem}}
\newcommand{\ethm}{\end {theorem}}
\newcommand{\bco}{\begin {coro}}
\newcommand{\eco}{\end {coro}}
\newcommand{\bex}{\begin {ex}}
\newcommand{\eex}{\end {ex}}

\newcommand{\be}{\begin {equation}}
\newcommand{\ee}{\end {equation}}
\newcommand{\bp}{\begin {proof}}
\newcommand{\ep}{\end {proof}}
\newcommand{\bee}{\begin {equation*}}
\newcommand{\eee}{\end {equation*}}
\newcommand{\rt}{\rightarrow}
\newcommand{\lb}{\label}
\newcommand{\ot}{\otimes}
\newcommand{\al}{\alpha}

\begin{document}
\title{Quaternionic slice regular functions and quaternionic Laplace transforms
\\[10mm]
}
\author{Gang Han\\School of Mathematics \\ Zhejiang University\\mathhgg@zju.edu.cn \\[3mm]
}
\date{June 12, 2010 }
\maketitle
\begin{abstract}
The functions studied in the paper are quaternion-valued functions of a quaternionic variable. It is show that the  left slice regular functions and right slice regular functions are related by a particular involution. The relation between  left slice regular functions, right slice regular functions and intrinsic regular functions is revealed.

The classical Laplace transform can be naturally generalized to quaternions in two different ways,  which transform a quaternion-valued function of a real variable to a left or right slice regular quaternion-valued function of a quaternionic variable.  The usual properties of the classical Laplace transforms are generalized to quaternionic Laplace transforms.

\end{abstract}

\textbf{ MAC2010}: 30G35, 35A22

\section{Introduction}
 \setcounter{equation}{0}\setcounter{theorem}{0}

The theory of slice regular functions (originally Cullen-regular functions) of a quaternionic variable was
born around 2005, and has rapidly developed since then. Most of the known results and applications of slice regular functions can be found in \cite{reg}. Many results in the theory of the complex analytic functions are generalized to slice regular functions, such as the Cauchy Integral Formula, the Maximum Modulus Principles, the Open Mapping Theorem, the Laurent series, and etc.

The theory of slice regular functions has now been applied to develop a new functional
calculus in a noncommutative setting \cite{css1}. The theory of regular functions has also been generalized to real alternative algebras \cite{gp}.

The set of left (resp. right) slice regular functions are not closed under the usual product of functions, but they are closed under the regular product ($*$-product). With respect to this regular product and the usual addition and scalar multiplication, the set of left (resp. right) slice regular functions (on some domain of $H$) becomes a real noncommutative associative algebra. The theory of left slice regular functions and that of right slice regular functions are parallel. Sometimes one may meet left slice regular functions and right slice regular functions at the same time. The Cauchy formulas for the left and for the right slice regular
functions involve different Cauchy kernels \cite{cgs}, unlike what happens in the classical case. The quaternionic Cauchy kernel

\[\phi(q,s)=-\left(|s|^{2}-2q \operatorname{Re}(s)+q^{2}\right)^{-1}(q-\bar{s})\] is a left slice regular function of $q$, as well as a right slice regular function of $s$ \cite{cgs}.

In the paper we will demonstrate the relation between left slice regular functions and right slice regular functions. They are related by an involution on the space of quaternionic functions.

Let $\mathbb{H}$ be the real algebra of quaternions and $U\subseteq \mathbb{H}$ be an axially symmetric slice domain and $C(U, \mathbb{H})$ be the set of continuous functions from $U$ to $\mathbb{H}$. Let $\mathcal{R}^0(U)$,  $\mathcal{R}^l(U)$ and $\mathcal{R}^r(U)$ be respectively the set of intrinsic regular functions , left regular functions and right regular functions on $U$. Let
\[ \eta: C(U, \mathbb{H})\rt C(U, \mathbb{H}), \eta(f)(q)=\overline{f(\overline{q})}.\] Then $\eta$ is an involution on $C(U, \mathbb{H})$, which maps left slice regular functions to right slice regular functions and vice versa, which induces an anti-isomorphism between the associative algebras $\mathcal{R}^l(U)$ and $\mathcal{R}^r(U)$.

Another main result is the isomorphism of associative algebras
\[\phi:\mathcal{R}^0(U)\bigotimes_\mathbb{R} \mathbb{H}\rt \mathcal{R}^l(U), \] through which some known results, including the slice derivatives and the regular conjugation of slice regular functions, can be better understood. See Theorem \ref{l} and Theorem \ref{a1}.

\bigskip

Next we introduce the quaternionic Laplace transform.

The classical Laplace transform is an integral transform named after its inventor Pierre-Simon Laplace. It transforms a complex-valued function of a real
variable to a complex-valued function of a complex variable. The transform has many applications in mathematics, physics and engineering.

The classical Laplace transform can be naturally generalized to quaternions, see \cite{l1} and \cite{l2}. As $\mathbb{H}$ is noncommutative, there are 2 types of quaternionic Laplace transforms (see Section 4),  which transform a quaternion-valued function of a real
variable to a left or right slice regular quaternion-valued function of a quaternionic variable.

The quaternionic Laplace transform can be generalized to quaternionic linear operators. The relation between quaternionic evolution operator and quaternionic Laplace transform were introduced in \cite{cs}.

In the paper, we show that the classical properties of the classical Laplace transforms can be generalized to quaternionic Laplace transforms. For example,  if $f(t)$ and $g(t)$  are  quaternion-valued functions defined for all real numbers $t\geq 0$, and $f \circ g$ is the convolution of $f$ and $g$. Then the left Laplace transform of $f \circ g$ is the regular product of the left Laplace transform of $f$ and the left Laplace transform of $g$ (Proposition \ref{m}).

Most of the notations in the paper are standard, and some notations will be explained when they first appear in the paper.

\bigskip

\textbf{Acknowledgements}\bigskip

I would like to heartily thank Hao Yu and Yeshun Sun for very helpful discussions with them.

\section{Some preliminary results on quaternions and slice regular functions}
 \setcounter{equation}{0}\setcounter{theorem}{0}

1. The quaternions\bigskip

The quaternions $\mathbb{H}$ were introduced by Hamilton in $1843,$ adding a multiplicative structure to $\mathbb{R}^{4} .$ It is a real division algebra. Let $\{1,i,j,k\}$ be the standard basis of $\mathbb{H}$.
 Let $\alpha:\mathbb{H}\rt \mathbb{H}, q\rt \overline{q}$ be the usual conjugation, which is an involution of $\mathbb{H}$. In the paper, an involution of an associative algebra means an anti-automorphism of the algebra of order two. Let $\mathbb{R}\cdot 1$ and $\operatorname{Im}(\mathbb{H})$ be respectively the $\mathbb{R}$-span of 1 and the $\mathbb{R}$-span of $\{i,j,k\}$, which are the respective  1 and -1 eigenspace of $\al$. One has $\operatorname{Im}(\mathbb{H})=\{q\in \mathbb{H}|q^2\in \mathbb{R}\ and\ q^2\leq 0\}$. One has the decomposition of $\mathbb{H}$: $$\mathbb{H}=\mathbb{R}\cdot 1\oplus \operatorname{Im}(\mathbb{H}).$$ For $q=a+bi+cj+dk\in \mathbb{H}$, let $\operatorname{Re}(q)=a$ and $\operatorname{Im}(q)=bi+cj+dk$, which are the respective real and imaginary part of $q$.
 Let $|q|^2=q\overline{q}$, which makes $\mathbb{H}$ into an 4-dimensional Euclidean space with the inner product $(p,q)=\operatorname{Re} (p\overline{q})$. The elements in $\mathbb{H}$ with unit norm constitute a group under multiplication, denoted by $Sp(1)$.

 Let $ \text{Aut}(\mathbb{H})$ be the automorphism group of $\mathbb{H}$.  One knows that  $\text{Aut}(\mathbb{H})\cong SO(3)$.

 Let  $\widetilde{\text{Aut}}(\mathbb{H})=\text{Aut}(\mathbb{H})\bigsqcup \al\cdot \text{Aut}(\mathbb{H})$ be the group generated by  $\text{Aut}(\mathbb{H})$ and $\al$, which consists of all the automorphisms and anti-automorphisms of $\mathbb{H}$. Since $\alpha$ restricts to minus identity on $\operatorname{Im} \mathbb{H}$, $\widetilde{\text{Aut}}(\mathbb{H})\cong O(3)$, and we will identify $\widetilde{\text{Aut}}(\mathbb{H})$ with $O(3)$.

Let $$\mathbb{S}=\{q\in \mathbb{H}|q^2=-1\},$$ which is the 2-dimensional unit sphere of purely imaginary elements in $\mathbb{H}$ with unit norm. It is the orbit of $i$ under $\text{Aut}(\mathbb{H})$. For any $I\in S$, let \[\mathbb{C}_I=\{x+Iy|x,y\in \mathbb{R}\},\] which is a subalgebra of $\mathbb{H}$ isomorphic to $\mathbb{C}$. For any $I,J\in \mathbb{H}$, if $J\neq \pm I$ then $\mathbb{C}_I\cap \mathbb{C}_J=\mathbb{R}$. One has that $\mathbb{C}_I=\mathbb{C}_J$ if and only if $J=\pm I$. One has $\mathbb{H}=\bigcup_{I\in S} \mathbb{C}_I$. This is called the slice structure of $\mathbb{H}$. For any open set $U$ in $\mathbb{H}$ and $I\in \mathbb{S}$, let \be\lb{y}U_I=U\cap \mathbb{C}_I.\ee

For any $I\in \mathbb{S}$, let $\tau_{I}:\mathbb{C} \rightarrow \mathbb{H}, x+iy\mapsto x+Iy$ be the embedding, whose image is $\mathbb{C}_I$. We also use $\tau_{I}$ to denote $\mathbb{C} \rightarrow \mathbb{C}_I, x+iy\mapsto x+Iy$.  $\tau_{I}$ is conjugation equivariant. As $\mathbb{C}_I=\mathbb{C}_{-I}$, there are 2 identifications $\mathbb{C}\rt \mathbb{C}_I$, $x+iy\mapsto x+Iy$ and $x+iy\mapsto x-Iy$.

The Galois group of $\mathbb{C}$ over $\mathbb{R}$ is the order 2 cyclic group generated by the complex conjugation $z\mapsto \overline{z}$, and will be  denoted by $C_2$.

One has the following homeomorphism of orbit spaces

$$
 \mathbb{H} / {SO}_{3} \cong \mathbb{C} / {C}_{2}.
$$

The map
$$\mathbb{C}\times \mathbb{S}\rt \mathbb{H}, (x+iy,J)\mapsto x+Jy,$$ is a surjective continuous map, whose restriction $(\mathbb{C}\setminus \mathbb{R })\times \mathbb{S}\rt \mathbb{H}\setminus \mathbb{R}$ is a double covering.\bigskip

 2. Definition of slice regular functions\bigskip

Recall the $U_I$ defined in (\ref{y}).

 \begin{defi} Let $U$ be an open set in $\mathbb{H}$ and $f: U \rightarrow \mathbb{H}$ be real differentiable. For every $I \in \mathbb{S},$ let $f_{I}:U_I\rt \mathbb{H}$ be the restriction of $f$ to $U_I$. The function $f$ is said to be left slice regular (left regular) if for every $I \in \mathbb{S},$
$$
\bar{\partial}_{I} f_{I}(x+I y):=\frac{1}{2}\left(\frac{\partial}{\partial x}+I \frac{\partial}{\partial y}\right) f_{I}(x+I y)=0
$$
on $U_I .$
Analogously, a function is said to be right slice regular (right regular) in $U$ if
$$
\bar{\partial}^r_{I} f_{I}(x+I y):=\frac{1}{2}\left(\frac{\partial}{\partial x} f_{I}(x+I y)+\frac{\partial}{\partial y} f_{I}(x+I y) I\right)=0
$$
on $U_I$.
\end{defi} The class of left (resp. right) regular functions on $U$ will be denoted by $\mathcal{R}^{l}(U)$ (resp. $\mathcal{R}^{r}(U)$).

The left (slice) $I$ -derivative of $f\in \mathcal{R}^{l}(U) $ at a point $q=x+I y$ is given by
$$
\partial_{I} f_{I}(x+I y)=\frac{1}{2}\left(\frac{\partial}{\partial x} f_{I}(x+I y)-I \frac{\partial}{\partial y} f_{I}(x+I y)\right)
$$
Similarly, the right (slice) $I$ -derivative of $f\in \mathcal{R}^{r}(U)$ at $q=x+I y$ is given by
$$
\partial_{I}^r f_{I}(x+I y)=\frac{1}{2}\left(\frac{\partial}{\partial x} f_{I}(x+I y)-\frac{\partial}{\partial y} f_{I}(x+I y) I\right)
$$
Let us now introduce the slice derivative for regular functions.
\begin{defi}  Let $U$ be an open set in $\mathbb{H}$. If $f\in \mathcal{R}^{l}(U)$, then the left slice derivative, $\partial f$, of $f,$ is defined by:
$$
\partial f(q)=\left\{\begin{array}{ll}
\partial_{I}f_I(q) & \text { if } q=x+I y, \quad y \neq 0 \\
\frac{\partial f}{\partial x}(x) & \text { if } q=x \in \mathbb{R}.
\end{array}\right.
$$
If $f\in \mathcal{R}^{r}(U)$, then the right slice derivative, $\partial^r f$, of $f,$ is defined by:
$$
\partial^rf(q)=\left\{\begin{array}{ll}
\partial_{I}^rf_I(q) & \text { if } q=x+I y, \quad y \neq 0 \\
\frac{\partial f}{\partial x}(x) & \text { if } q=x \in \mathbb{R}.
\end{array}\right.
$$
\end{defi}

Note that if $f$ is a left regular function, then its derivative is still left regular because
\[
\bar{\partial}_{I}\left(\partial f(x+I y)\right)=\partial\left(\bar{\partial}_{I} f(x+I y)\right)=0.
\]

 Similar result holds for right regular functions.

 Let $r>0$ and $B(0, r)$ be the open ball centered at $0$ with radius $r$. A function $f : B(0, r)\rt \mathbb{H}$ is left (resp. right) regular on $B(0, r)$ if and only if $f$ has
a series representation of the form $f(q)=\sum_{n=0}^\infty q^n a_n$ (resp. $f(q)=\sum_{n=0}^\infty a_n q^n $), where
$\{a_n\}, n \in \mathbb{N},$ is a sequence of quaternions satisfying $(limsup_{n\rt \infty} |a_n|^{1/n})^{-1}\geq r$. See Theorem 2.1.5.Theorem 2.1.6. \cite{ent}.\bigskip

3. Intrinsic regular functions \bigskip

 Let $D\subseteq \mathbb{C}$ be a domain symmetric with respect to the real axis. Then $D\cap \mathbb{R}\neq \emptyset$.  Let $f:D\mapsto \mathbb{C}$ be a complex function.  Then  $f$ is called intrinsic if $f$ is $C_2$-equivariant, i.e.,
 \be\lb{a}
f(\bar{z})=\overline{f(z)}.
\ee

 \begin{defi}  Let $U \subseteq \mathbb{H}$. We say that $U$ is axially symmetric if $U$ is $SO(3)$-invariant, i.e., for every $x+I y \in U,$ all the elements $x+\mathbb{S} y=\{x+J y | J \in \mathbb{S}\}$ are contained in U. We say that $U$ is a slice domain if it is a connected open set whose intersection with every complex plane $\mathbb{C}_{I}$ is connected. \end{defi}

 If $U\subseteq \mathbb{H}$ is an axially symmetric slice domain, then $U\cap \mathbb{R}\neq \o$.

Let $U\subseteq \mathbb{H}$ be an axially symmetric slice domain and $C(U, \mathbb{H})$ be the set of continuous functions from $U$ to $\mathbb{H}$. Let
\be\lb{h} \eta: C(U, \mathbb{H})\rt C(U, \mathbb{H}), \eta(f)(q)=\overline{f(\overline{q})}.\ee

Then one has \[\eta^2=1\] and \[\eta(f\cdot g)=\eta(g)\cdot\eta(f),\ \ f,g\in C(U, \mathbb{H}),\] where $f\cdot g$ is the pointwise multiplication. So $\eta$ is an involution on $C(U, \mathbb{H})$. This involution $\eta$ plays an important role in the paper!

 \begin{defi} Let $U$ be an axially symmetric slice domain. A function $F:U\rt \mathbb{H}$ is called \textit{slice-preserving}, if $F(U\cap \mathbb{C}_{I})\subseteq U\cap \mathbb{C}_{I}$ for each $I\in \mathbb{S}$. If $F$ is slice-preserving and left regular, then $F$ is called \textit{intrinsic regular}. The set of intrinsic regular functions on $U$ will be denoted by $\mathcal{R}^0(U)$. \end{defi}

The following result is clear.

\begin{prop}
Assume that $F:U\rt \mathbb{H}$ is  slice-preserving, where $U$ is an axially symmetric slice domain in $\mathbb{H}$, then $F$ is left regular if and only if $F$ is right regular if and only if $F$ is intrinsic regular.
\end{prop}

By Remark $1.31$ of \cite{reg}, one know that if  $f$ is a left regular function on $B(0, r)$ (for some $r>0$ ), then $f$ is slice preserving if, and only if, the power series expansion $f(q)=\sum_{n \in \mathbb{N}} q^{n} a_{n}$ has real coefficients $a_{n} \in \mathbb{R}$.

If $f$ is an intrinsic regular function, then $\partial(f)(q)= \partial^r(f)(q)$, which will be denoted by $f^{'}(q)$, i.e., for intrinsic regular functions $f$, \[f^{'}(q)=\partial(f)(q)= \partial^r(f)(q).\]

\bigskip

4. Extension from complex intrinsic holomorphic functions to quaternionic intrinsic regular functions\bigskip

The following construction is due to Feuter \cite{f}.

 Let $f(x+iy)=u(x,y)+iv(x,y)$ be an intrinsic complex function on  a domain $D\subseteq \mathbb{C}$ symmetric with respect to the real axis, where $u$ and $v$ are real-valued functions. Let $\widetilde{D}\subseteq \mathbb{H}$ be the union of $SO(3)$-orbit of all the points in $D$, i.e., $\widetilde{D}= \{x+\mathbb{S} y| x+iy\in D\}$. Then $f$ induces a function, $\widetilde{f}$, from $\widetilde{D}$ to $\mathbb{H}$, such that for any $J\in \mathbb{S}$,

\be\lb{a}
\widetilde{f}:\widetilde{D}\rt \mathbb{H}, \widetilde{f}(a+J b)=u(a, b)+J v(a, b).
\ee
 This is well-defined because $f$ is intrinsic. This induced function $\widetilde{f}$ is $O(3)$-equivariant and slice-preserving. This construction is later generalized to regular functions.

Assume further that $f$ is an intrinsic holomorphic function on $D$, then $\widetilde{f}$ is left and right regular on $\widetilde{D}$ preserving each $\widetilde{D}\cap \mathbb{C}_{I}$, thus is an {intrinsic regular function} on $\widetilde{D}$.  Let $ H(D, \mathbb{C})$ be the set of intrinsic holomorphic function on $D$ and define

\be\lb{i} \text{ext}: H(D, \mathbb{C})\rt \mathcal{R}^0(\widetilde{D}), f\mapsto \widetilde{f},  \ee which is an isomorphism of commutative and associative algebras.
Let $U=\widetilde{D}$. Then the inverse of this map is the restriction of $g\in  \mathcal{R}^0(U)$ to some single slice $U_I=U\cap \mathbb{C}_I$ (for some $I\in \mathbb{S}$), where $g|_{U_I}\in H(U_I, \mathbb{C}_I)$. Write $H(U_I)=H(U_I, \mathbb{C}_I)$. Then (\ref{i}) can also be written as \be\lb{z} \text{ext}: H(U_I)\rt  \mathcal{R}^0(U). \ee

\bex
The exponential function
$$
\exp : \mathbb{H} \rightarrow \mathbb{H}, \quad q \rightarrow e^q= \sum_{n=0}^{\infty} \frac{q^{n}}{n !}
$$

\noindent is an intrinsic regular function induced from the usual complex exponential function.
\eex

\section{Main results on regular functions}
 \setcounter{equation}{0}\setcounter{theorem}{0}

In this section $U$ will always be an axially symmetric slice domain in $\mathbb{H}$.

\begin{prop}\lb{b}(Proposition 2.2.7 and Corollary 2.2.8 of \cite{ent}) Let $f \in \mathcal{R}^l(U).$ For any orthogonal unit vectors $I, J, K \in \mathbb{S}$ with $IJ=K$, there exist 4 uniquely determined
intrinsic regular function $h_{m}$ on $U$, such that
$ f(q)=h_{0}(q)+h_{1}(q) I+h_{2}(q) J+h_{3}(q) K
$, $q\in U$. Each $h_{m}$ is a real analytic function.\end{prop}

The  following result is a corollary of the Identity Principle (Theorem 1.12 of \cite{reg}).
\begin{prop}\lb{g}
Assume $f, g\in R
^l(U)$ and $W\subseteq U\cap \mathbb{R}$ is a nonempty open subset. If $f(q)=g(q)$ for any $q\in W$, then $f=g$.
\end{prop}

For any $f\in \mathcal{R}^l(U)$, one knows that $f$ can be recovered from its restriction to some single slice $U_I=U\cap \mathbb{C}_I$ with $I\in \mathbb{S}$, see Corollary 1.16 of \cite{reg}. In fact, more is true.

Let $U$ be an  axially symmetric slice domain of $\mathbb{H}$ and $f\in \mathcal{R}^l(U)$. Then $f$ can be recovered from its restriction to $U\cap \mathbb{R}$ as follows.
By Proposition \ref{b}, one knows that $f(q)=f_{0}(q)+f_{1}(q) i+f_{2}(q) j+f_{3}(q) k$ with $f_m\in \mathcal{R}^0(U)$. Fix $I\in \mathbb{S}$. Each $f_m(q)$ is determined by its restriction, $\bar{f_m}$, to $U_I$. As $\bar{f_m}$ is holomorphic on $U_I$, $\bar{f_m}$  is determined by its restriction to $U\cap \mathbb{R}$, which is denoted by $f_m|_\mathbb{R}$. In fact, $\bar{f_m}$  is the analytic continuation of  $f_m|_\mathbb{R}$ to
$U_I$, and $f_m=\text{ext}(\bar{f_m})$ (see (\ref{a})).

\begin{prop}\lb{j}
 A left (resp. right) regular function $f:U\rt \mathbb{H}$ is intrinsic regular if and only if $f(U\cap \mathbb{R})\subseteq \mathbb{R}$.
\end{prop}
\bp
We will prove the result in the case $f\in \mathcal{R}^l(U)$.

If $f$ is intrinsic regular, then $f(U\cap \mathbb{C}_{I})\subseteq U\cap \mathbb{C}_{I}$ for each $I\in \mathbb{S}$. Since $\mathbb{C}_{I}\cap \mathbb{C}_{J}=\mathbb{R}$ for $J\neq \pm I$, $f(U\cap \mathbb{R})\subseteq \mathbb{R}$.

Conversely, assume $f(U\cap \mathbb{R})\subseteq \mathbb{R}$ and
$ f(q)=f_{0}(q)+f_{1}(q) I+f_{2}(q) J+f_{3}(q) K
$ with each $f_{m}\in \mathcal{R}^0(U)$. Then $f_{m}(q)=0$ for $q\in U\cap \mathbb{R}$ for $m=1,2,3$. Since $f_{m}(q)$ is complex holomorphic on each slice $U_{I}$, $f_{m}(q)=0$ for $q\in U_{I}$ for  $m=1,2,3$. As $U=\cup_{I\in \mathbb{S}} U_I$, $f_{m}(q)=0$ for each $q\in U$ for  $m=1,2,3$, and $f=f_0$ is intrinsic regular.

\ep
\bigskip

For any $f,g\in \mathcal{R}^l(U)$, the function $h(q)=f(q)g(q)$ may not be in $\mathcal{R}^l(U)$. The left regular product $f*g$ of $f$ and $g$ (see Definition 1.27 of \cite{reg}), is also in  $\mathcal{R}^l(U)$.   If $f(q)=\sum_{n \in \mathbb{N}} q^{n} a_{n}$ and $g(q)=\sum_{n \in \mathbb{N}} q^{n} b_{n}$ are left regular functions on $B(0, r)$. Then the regular product of $f$ and $g$  is the regular function defined by
\be \lb{f}
(f * g)(q)=\sum_{n \in \mathbb{N}} q^{n} \sum_{k=0}^{n} a_{k} b_{n-k}
\ee
on the same ball $B(0, r)$.

Assume $f\in \mathcal{R}^l(U)$. For the definition of its regular conjugate $f^{c}$ and its symmetrization $f^s$, see Definition 1.33 of \cite{reg}.

  Assume $f\in \mathcal{R}^l(U)$ and $f \neq 0$. The regular reciprocal $f^{-*}$ of $f$ is the function defined on $U \backslash Z_{f^{s}}$, where $Z_{f^{s}}=\{q\in \mathbb{H}|f^s(q)=0\}$, by
\[
f^{-*}=\frac{1}{f^{s}} f^{c}.
\]
See Definition 5.1 of \cite{reg}.

\begin{prop}\lb{k}
Assume $f,g\in \mathcal{R}^l(U)$. One has

(1) The left regular product $f*g$ is the left regular function in $\mathcal{R}^l(U)$ such that $(f*g)(q)=f(q)g(q)$ for $q\in U\cap \mathbb{R}$.

(2) The regular conjugate $f^{c}$ of $f$ is the function in $\mathcal{R}^l(U)$ such that for $q\in U\cap \mathbb{R}$, $f^{c}(q)=\overline{f(q)}$.

(3) The symmetrization $f^s$ of $f$ is the intrinsic regular function in $\mathcal{R}^0(U)$ such that for $q\in U\cap \mathbb{R}$, $f^{s}(q)=f(q)\overline{f(q)}$.

(4) The regular reciprocal $f^{-*}$ of $f$ is the function in $\mathcal{R}^l(U)$ such that for $q\in (U\cap \mathbb{R})\backslash Z_{f^s}$, $f^{-*}(q)=f(q)^{-1}$.
\end{prop}
\bp

By Proposition \ref{g}, we only need to show that the value of the regular functions $f*g$, $f^{c}$, $f^s$ and $f^{-*}$ in $q\in U\cap \mathbb{R}$ is as stated in the proposition.

(1) Without loss of generality one can assume $U\supseteq B(0, r)$ for some $r>0$. Then by (\ref{f}), $(f*g)(q)=f(q)g(q)$ for $q\in B(0, r)\cap \mathbb{R}$. By Proposition \ref{b}, the function $ U\cap \mathbb{R}\rt \mathbb{H}, q\mapsto  f(q)g(q)$ is a  function of a real variable with each component real analytic, thus there is a unique left regular function $h$ such that $h(q)=f(q)g(q)$ for $q\in U\cap \mathbb{R}$. As $h(q)=(f*g)(q)$ for $q\in B(0, r)\cap \mathbb{R}$, $h=f*g$ on $U$  by Proposition \ref{g}.

(2) It can be proved similarly as in (1).

(3) It follows from (1) and (2) as $f^s=f*f^c$.

(4) For $q\in (U\cap \mathbb{R})\backslash Z_{f^s}$, $f(q)\cdot \frac{1}{f^{s}(q)} f^{c}(q)=1$, thus $f^{-*}(q)=\frac{1}{f^{s}(q)} f^{c}(q)=f(q)^{-1}$. Then (4) follows from Proposition \ref{g}.

\ep

\bigskip

Now let us review the Splitting Lemma before we prove the next result. Let $f:U\rt \mathbb{H}$ be a real differentiable function with $U$ an axially symmetric slice domain.  Then $f$ is left regular if and only if for any $I \in \mathbb{S}$ and any $J \in \mathbb{S}$ with $J \perp I,$ there exist two holomorphic functions
$F, G: U_{I} \rightarrow \mathbb{C}_{I}$ such that for every $z=x+y I,$
$$
f_{I}(z)=F(z)+G(z) J.
$$
And, $f$ is  right regular if and only if for any $I \in \mathbb{S}$ and any $J \in \mathbb{S}$ with $J \perp I,$ there exist two holomorphic functions
$F, G: U_{I} \rightarrow \mathbb{C}_{I}$ such that for every $z=x+y I,$ $$f_{I}(z)=F(z)+ J G(z).$$

The space $\mathcal{R}^l(U)$ (resp. $\mathcal{R}^r(U)$) is a real associative algebra equipped with the regular product. Recall the involution $\eta$ defined in (\ref{h}), let \[C(U,\mathbb{H})^\eta=\{f\in C(U,\mathbb{H})|\eta(f)=f \}.\]
\begin{theorem}\lb{l}
Let $U$ be an axially symmetric slice domain of $\mathbb{H}$.

(1) If $f\in \mathcal{R}^l(U)$ then $\eta(f)\in \mathcal{R}^r(U)$ and viceversa.  One has $\mathcal{R}^0(U)= \mathcal{R}^l(U)\cap C(U,\mathbb{H})^\eta= \mathcal{R}^r(U)\cap C(U,\mathbb{H})^\eta$.

(2) The map \[\phi:\mathcal{R}^0(U)\bigotimes_\mathbb{R} \mathbb{H}\rt \mathcal{R}^l(U), f(q)\otimes \lambda\mapsto f(q) \lambda, \] is an isomorphism of real associative algebras, where the product of $f_0\otimes I$ and $g_0\otimes J$ in $\mathcal{R}^0(U)\bigotimes_\mathbb{R} \mathbb{H}$ is $ f_0\otimes I \cdot g_0\otimes J=f_0  g_0\otimes IJ$. One has that $\mathcal{R}^0(U)$ is the center of the associative algebra $\mathcal{R}^l(U)$.

The map $\phi$ is also an isomorphism of a left (free) $\mathcal{R}^0(U)$-module and a right $\mathbb{H}$-module. 

The results for $\mathcal{R}^r(U)$ are symmetric.
\end{theorem}

 \bp
 (1)Let $f(q)\in \mathcal{R}^l(U)$. By the Splitting Lemma, for any $I \in \mathbb{S}$ and any $J \in \mathbb{S}$ with $J \perp I,$ there exist two holomorphic functions
$F, G: \Omega_{I} \rightarrow L_{I}$ such that for every $z=x+y I,$ $f_{I}(z)=F(z)+G(z) J$.

Let $g(q)=\overline{f(\overline{q})}$, then for every $z=x+y I,$ $g_I(z)=\overline{f_I(\overline{z})}=\overline{F(\overline{z})+G(\overline{z})J}=\overline{F(\overline{z})}+J(-\overline{G(\overline{z})})$. Since
$F, G: \Omega_{I} \rightarrow L_{I}$ are holomorphic functions, $\overline{F(\overline{z})}$ and $-\overline{G(\overline{z})}$ are also holomorphic, thus $g(q)=\overline{f(\overline{q})}$ is right regular.

One has $\mathcal{R}^0(U)\subseteq \mathcal{R}^l(U)\cap C(U,\mathbb{H})^\eta$ by definition. Conversely, if $f\in  C(U,\mathbb{H})^\eta$ then for $q\in U\cap \mathbb{R}$, $ \overline{f(\overline{q})}=\overline{f(q)}=f(q)$ thus $f(q)$ is real. So, if $f\in \mathcal{R}^l(U)\cap C(U,\mathbb{H})^\eta$, then $f\in \mathcal{R}^0(U)$ by Proposition \ref{j}. So $\mathcal{R}^0(U)\supseteq \mathcal{R}^l(U)\cap C(U,\mathbb{H})^\eta$.

 (2) Proposition \ref{b} implies that
\[
\mathcal{R}^l(U)=\mathcal{R}^0(U) \oplus \mathcal{R}^0(U) i \oplus \mathcal{R}^0 j\oplus \mathcal{R}^0(U) k,
\]

 which implies that $\phi$ is an $\mathbb{R}$-linear isomorphism of $\mathcal{R}^0(U)\otimes_\mathbb{R} \mathbb{H}$ and $\mathcal{R}^l(U)$. It is directly verified that $\phi$ preserves multiplications in the two associative algebras. So $\phi$ is  an isomorphism of real associative algebras. It is easy to see that $\mathcal{R}^0(U)$ is the center of the associative algebra $\mathcal{R}^l(U)$.

 The left $\mathcal{R}^0(U)$-module structure and right $\mathbb{H}$-module structure on $\mathcal{R}^l(U)$  and $\mathcal{R}^0(U)\otimes_\mathbb{R}\mathbb{H}$ are the obvious ones. It is  clear that $\phi$ is an isomorphism of the left $\mathcal{R}^0(U)$-module and the right $\mathbb{H}$-module.

\ep

{In the rest, we will identify $\mathcal{R}^l(U)$ with $\mathcal{R}^0(U)\bigotimes_\mathbb{R} \mathbb{H}$ by $\phi$}.

\begin{theorem}\lb{a1} Let $U$ be an axially symmetric slice domain of $\mathbb{H}$.

 (1)One has $\eta(f*g)=\eta(g)*\eta(f)$ for $f,g\in \mathcal{R}^l(U)$, where the 1st (resp. the 2nd) product is left (resp. right) regular product. So $\eta:\mathcal{R}^l(U)\rt \mathcal{R}^r(U)$ is an anti-isomorphism of associative algebras. Similarly, $\eta:\mathcal{R}^r(U)\rt \mathcal{R}^l(U)$ is also an anti-isomorphism of associative algebras.

(2) The map  \[\mathcal{R}^0(U)\bigotimes_\mathbb{R} \mathbb{H}\rt \mathcal{R}^0(U)\bigotimes_\mathbb{R} \mathbb{H}, f\otimes \mu\mapsto f\otimes \overline{\mu}\] induces the regular conjugation on $\mathcal{R}^l(U)$, i.e., $(f\otimes \mu)^c=f\otimes \overline{\mu}$. And similar result holds for $\mathcal{R}^r(U)$.

(3)Assume $f\in \mathcal{R}^l(U)$. Write $f=\sum_{m=0}^3 f_m(q)\otimes J_m$, with $\{J_m| m=0,1,2,3\}$ a basis of $H$. Then \[\partial(\sum_{m=0}^3 f_m(q)\otimes J_m)=\sum_{m=0}^3 f_m^{'}(q)\otimes J_m.\]

(4) Let $I\in \mathbb{S}$ and recall the map $\text{ext}$ in (\ref{z}). Let $\psi: H(U_I)\bigotimes_\mathbb{R} \mathbb{H}\rt \mathcal{R}^l(U)$ be the composition of \[\text{ext}\otimes 1:H(U_I)\bigotimes_\mathbb{R} \mathbb{H}\rt \mathcal{R}^0(U)\bigotimes_\mathbb{R} \mathbb{H}\] with $\phi:\mathcal{R}^0(U)\bigotimes_\mathbb{R} \mathbb{H}\rt \mathcal{R}^l(U)$. Then $\psi$ is an isomorphism of associative algebras.

If one identifies any element  $g\otimes J\in H(U_I)\bigotimes_\mathbb{R} \mathbb{H}$ with the function $U_I\rt H, q\mapsto g(q)J$, then $\psi(g\ot J)$ is the regular extension of $g\ot J$. (For the definition of regular extension, see Lemma 1.22 of \cite{reg}.) Similar result holds for $\mathcal{R}^r(U)$.
 \end{theorem}

\bp

(1)Assume $f,g\in \mathcal{R}^l(U)$ with $f(q)=f_0(q)\otimes I$ and $g(q)=g_0(q)\otimes J$.
Then
 $$
\begin{aligned}
\eta(f*g)&=\eta(f_0\otimes I \cdot g_0\otimes J  )\\
&=\eta(f_0 g_0\otimes IJ )\\
&=\overline{IJ}\otimes f_0 g_0\\
&=\overline{J}\otimes g_0 \cdot \overline{I}\otimes f_0\\
&=\eta(g)*\eta(f).
\end{aligned}
$$
Then it follows by linearity that   $\eta(f*g)=\eta(g)*\eta(f)$ for any $f,g\in \mathcal{R}^l(U)$. Similarly one also has $\eta(f*g)=\eta(g)*\eta(f)$ for any $f,g\in \mathcal{R}^r(U)$. Since $\eta^2=1$, $\eta:\mathcal{R}^l(U)\rt \mathcal{R}^r(U)$ is an anti-isomorphism with inverse $\eta:\mathcal{R}^r(U)\rt \mathcal{R}^l(U)$, which is also an anti-isomorphism.

(2) and (3) can be verified directly by definition.

(4) As $\phi$ and $\text{ext}\ot 1$ are both isomorphisms of associative algebras, their composition $\psi$ is also an isomorphism of associative algebras. Let $g\ot J\in H(U_I)\bigotimes_\mathbb{R} \mathbb{H}$, $\widetilde{g}=\text{ext}(g)\in \mathcal{R}^0(U)$ and $f=\psi(g\ot J)\in \mathcal{R}^l(U)$. Then for $q\in U$, $f(q)=\psi(g\ot J)(q)=\phi(\widetilde{g}\ot J)(q)=\widetilde{g}(q)J$. So for $q\in U_I$, \[f(q)=\widetilde{g}(q)J=g(q)J=(g\ot J)(q).\] Thus $f|_{U_I}=g\ot J$ and $f$ is the unique regular extension of $g\ot J$.

 \ep

\section{Left and right quaternionic Laplace transforms }
 \setcounter{equation}{0}\setcounter{theorem}{0}
Let $g(t)$ be a complex-valued function, defined for all real numbers $t\geq 0$. Then its classical Laplace transform $L\{g\}$ is defined by  \[L\{g\}(s)=\int_0^\infty e^{-ts}g(t)dt,\] for those complex numbers $s\in \mathbb{C}$ such that the integral convergents.  It will be referred as the complex Laplace transform in the paper. For the main properties of complex Laplace transform, please refer to the classical book \cite{it}.

Let $f(t)$  be a quaternion-valued function defined for all real numbers $t\geq 0$. Note that in last section $f(t)$ is a quaternion-valued function of a quaternionic variable, but in this section $f(t)$ is a quaternion-valued function of a \textbf{real} variable!

Define $$L^l\{f\}(s)=\int_0^\infty e^{-ts}f(t)dt,$$  for those quaternions $s\in \mathbb{H}$ such that the integral convergents,  to be the left quaternionic Laplace transform of $f$; and $$L^r\{f\}(s)=\int_0^\infty f(t)e^{-ts}dt,$$  for those quaternions $s\in \mathbb{H}$ such that the integral convergents, to be the right quaternionic Laplace transform of $f$.

Note that if $f$ is real-valued, then $L^r\{f\}=L^l\{f\}$.

Let $f(t)$  be a quaternion-valued function defined for all real numbers $t\geq 0$. Assume $a\geq 0$. The function $f(t)$  is said to be of exponential order $a$ on $t\geq 0$ if there exist two positive constants $K$ and $T$, such that for all $t>T$
$$
|f(t)| \leq K e^{a t}.
$$
Once we say that  $f(t)$  is of exponential order $a$ then we always assume that $a\geq 0$.

\begin{theorem} Let $f(t)$  be a quaternion-valued function defined for all real numbers $t\geq 0$, which is of exponential order $a$ (with $a\geq 0$), and is continuous or piecewise continuous in every finite interval $(0, T).$ Then

(1) The left Laplace transform $L^l\{f\}(s)$ of $f(t)$ convergent absolutely for all $s\in \mathbb{H}$ provided $\operatorname{Re}(s)>a$.

(2) The left Laplace transform
$L^l\{f\}(s)=\int_{0}^{\infty} e^{-st}f(t) d t$ is uniformly convergent with respect to $s\in \mathbb{H}$ provided $\operatorname{Re}(s) \geq a_{1}$ where $a_{1}>a$.

(3) The map $f\mapsto L^l\{f\}$, where $f$ is of exponential order $a,$  is right $\mathbb{H}$-linear.

We will denote $F=L^l\{f\}$ in the rest of the theorem.

(4)If $f(t)$ is real-valued, then $F(s)$ is intrinsic regular on $\operatorname{Re}(s)>a$  and  $F{'}(s)=L^l\{-tf(t)\}$. The restriction of $F$ to each $\mathbb{C}_I$ is just the complex Laplace transform of $f$.

(5)In general, $F(s)$ is left regular on  $\operatorname{Re} (s)>a$, and  $\partial F(s)=L^l\{-tf(t)\}$.

\end{theorem}

\bp
(1) and (2) can be proved as in the complex case.

(3) follows by definition.

(4) Let $U=\{s\in \mathbb{H}|\operatorname{Re} (s)>a\}.$  For any $I\in \mathbb{S}$, $F({U_I})\subseteq U_I$ and $F|_{U_I}$ is just the usual complex Laplace transform of $f$, which is holomorphic on $U_I$. So $F$ is intrinsic regular. For each $I\in \mathbb{S}$ and  $s\in U_I$ with $\operatorname{Re} (s)>a$, one has that $F{'}(s)=\int_{0}^{\infty} e^{-st}(-tf(t)) d t$ , thus for any $s\in U$, $F{'}(s)=L^l\{-tf(t)\}$.

(5) In general, $f=f_0+f_1 i+f_2 j+f_3 k$ where $f_m$ are real valued functions. Then by (3), $L^l\{f\}=L^l\{f_0\}+L^l\{f_1\} i+L^l\{f_2 \}j+L^l\{f_3\} k$, which is left regular as each $L^l\{f_m\}$ is intrinsic regular.
\ep
\begin{rem}
The results for right Laplace transform are completely symmetric and are omitted.
\end{rem}

For a quaternion-valued function $f$ defined for all real numbers $t\geq 0$, $\overline{f}$ is the function defined by $\overline{f}(t)=\overline{f(t)}$. Then $f$ is real-valued if and only if $f=\overline{f}$.
\begin{prop}Let $f(t)$ be a quaternion-valued function defined for all real numbers $t\geq 0$. Then
 $\eta(L^l\{f\})=L^r{\{\overline{f}\}}$ and  $\eta(L^r\{f\})=L^l{\{\overline{f}\}}$.
 \end{prop}
\bp  Write $f(t)=\sum_{m=0}^3 f_m(t)\cdot J_m$, where $J_0=1, J_1=i, J_2=j, J_3=k$ and $f_m(t)$ are all real-valued functions. Then $\overline{f}(t)=\sum_{m=0}^3 \overline{J_m} \cdot f_m(t)$.

 One has
 $$
\begin{aligned}
L^l\{f\}(s)&=\int_0^\infty e^{-ts}f(t)dt=\int_0^\infty e^{-ts}(\sum_{m=0}^3 f_m(t)\cdot J_m)dt\\
&=\sum_{m=0}^3 (\int_0^\infty e^{-ts}f_m(t) dt )\cdot  J_m=\sum_{m=0}^3 L^l\{f_m\}(s)\cdot  J_m
\end{aligned}
$$

 $$
\begin{aligned}
L^r\{\overline{f}\}(s)&=\int_0^\infty \overline{f}(t)e^{-ts}dt=\int_0^\infty(\sum_{m=0}^3 \overline{J_m}\cdot  f_m(t))e^{-ts}dt\\
&=\sum_{m=0}^3 \overline{J_m}\cdot  (\int_0^\infty e^{-ts}f_m(t) dt )=\sum_{m=0}^3 \overline{J_m}\cdot  L^l\{f_m\}(s)
\end{aligned}
$$

So  $\eta(L^l\{f\})=L^r\{\overline{f}\}$. The result  $\eta(L^r\{f\})=L^l\{\overline{f}\}$ can be similarly proved.

\ep

\bex

 Let $b\in \mathbb{H}$ and $f(t)=e^{bt}$. Let $F=L^l\{f\}$ and $G=L^r\{f\}$.

  For $s\in \mathbb{R}$ and $\operatorname{Re}(s)>\operatorname{Re}(b)$, $F(s)=G(s)=\int_0^\infty e^{-ts}e^{bt}dt=(s-b)^{-1}$.

   As $F$ is left regular, it is the left regular reciprocal of $s-b$. So $F(s)=(s^2-2\operatorname{Re}(b) s+ |b|^2)^{-1}(s-\overline{b})$, where $\operatorname{Re}(s)>\operatorname{Re}(b)$.

 As $G(s)$ is right regular, it is the right regular reciprocal of $s-b$. Then $G(s)=(s-\overline{b})(s^2-2\operatorname{Re}(b)s+ |b|^2)^{-1}$, where $\operatorname{Re}(s)>\operatorname{Re}(b)$.
\eex
In the rest of the paper we will consider only left Laplace transforms, and $f(t)$ will always be a quaternion-valued function defined for all real numbers $t\geq 0$.
\begin{prop}
 If ${L}\{f(t)\}=F(s),$ then
$$
\mathcal{L}^l\left\{e^{-a t} f(t)\right\}=F(s+a),
$$
where $a$ is a real constant.
\end{prop}
\bp We have, by definition,
$$
\mathcal{L}^l\left\{e^{-a t} f(t)\right\}=\int_{0}^{\infty} e^{-(s+a) t} f(t) d t=F(s+a).
$$
\ep
Note that one does not have $\mathcal{L}^l\left\{f(t)e^{-a t} \right\}=F(s+a)$ because of the noncommutativity of $\mathbb{H}$.

Recall the Heaviside function $$ H(x)=\left\{\begin{array}{ll}
0, & x<0; \\
1, & x\geq 0.
\end{array}\right.$$

\begin{prop} If $\mathcal{L}^l\{f(t)\}=F(s),$ then for $a>0$ one has:
$$
\mathcal{L}^l\{f(t-a) H(t-a)\}=e^{-a s} F(s)=e^{-a s} \mathcal{L}^l\{f(t)\}.
$$
\end{prop}
The proof is as in the complex case and is omitted.

\begin{prop}\lb{x}(Laplace Transforms of Derivatives). Let $f(t)$  be a quaternion-valued function defined for all real numbers $t\geq 0$, which is of exponential order $a$. Let $U=\{s\in \mathbb{H}| \operatorname{Re}(s)>a\}$. If $\mathcal{L}^l\{f(t)\}=F(s),$ then
$\mathcal{L}^l\left\{f^{\prime}(t)\right\}=s \mathcal{L}^l\{f(t)\}-f(0+)=s F(s)-f(0+)$ for those $s\in U$;
More generally one has, \[\mathcal{L}^l\left\{f^{(n)}(t)\right\}=s^{n} F(s)-s^{n-1} f(0+)-s^{n-2} f^{\prime}(0+)-\cdots-s f^{(n-2)}(0+)-f^{(n-1)}(0+)\]  for those $s\in U$, where $f^{(r)}(0+)=lim_{t\rt 0+}f^{(r)}(t)$, $r=0,1,\cdots, n-1$.
\end{prop}
\bp  Integrate by parts, one has
$$
\begin{array}{l}
\mathcal{L}^l\left\{f^{\prime}(t)\right\}=\int_{0}^{\infty} e^{-s t} f^{\prime}(t) d t\\
=\left[e^{-s t} f(t)\right]_{0}^{\infty}-\int_{0}^{\infty} [d(e^{-s t})\cdot f(t)] \\
=\left[e^{-s t} f(t)\right]_{0}^{\infty}+s \int_{0}^{\infty} e^{-s t} f(t) d t \\
=-f(0+)+s F(s),
\end{array}
$$
which holds for those $s\in \mathbb{H}$ with $\operatorname{Re}(s)>a$. The formula for $\mathcal{L}^l\left\{f^{(n)}(t)\right\}$ can be proved similarly.\ep

 Let $f(t)$ and $g(t)$  be  quaternion-valued functions defined for all real numbers $t\geq 0$. Let $f \circ g$ be the convolution of $f(t)$ and $g(t)$, which is defined by the integral
$$
(f\circ g)(t)=\int_{0}^{t} f(t-\tau) g(\tau) d \tau.
$$
Note that $(f\circ g)(t)$ may not equal to $(g\circ f)(t)$.

\begin{prop}\lb{m} (Convolution Theorem). Assume that $f(t)$ and $g(t)$ are of exponential order $a$ and $b$ respectively. Assume $\mathcal{L}^l\{f(t)\}=F(s)$ and $\mathcal{L}^l\{g(t)\}=G(s).$ Let $c=max\{a,b\}$ and $U=\{s\in \mathbb{H}| \operatorname{Re}(s)>c\}$. Then  $\mathcal{L}^l\{f\circ g\}\in \mathcal{R}^l(U)$ and $\mathcal{L}^l\{f\circ g\}(s)=(F*G)(s)$ for $s\in U$. In particular, if $s\in \mathbb{R}$ and $s>c$,
$$
\mathcal{L}^l\{f\circ g\}(s)=\mathcal{L}^l\{f(t)\} \mathcal{L}^l\{g(t)\}=F(s) G(s).
$$

\end{prop}
\bp
It is directly verified that for any $\varepsilon>0$, $(f\circ g)(t)$ is of exponential order $c+\varepsilon$. Then $\mathcal{L}^l\{f\circ g\}$ is defined on $U=\{s\in \mathbb{H}| \operatorname{Re}(s)>c\}$, and it is left regular on $U$.

By definition, for $s\in \mathbb{R}\cap U$,
$$
\begin{aligned}
F(s) G(s) &=\int_{0}^{\infty} e^{-s \sigma} f(\sigma) d \sigma \int_{0}^{\infty} e^{-s \mu} g(\mu) d \mu \\
&=\int_{0}^{\infty} \int_{0}^{\infty} e^{-s(\sigma+\mu)} f(\sigma) g(\mu) d \sigma d \mu.
\end{aligned}
$$

We make the change of variables $\mu=\tau, \sigma=t-\mu=t-\tau$. Consequently,  for $s\in \mathbb{R}\cap U$, it becomes
$$
\begin{aligned}
F(s) G(s) &=\int_{0}^{\infty} e^{-s t} d t \int_{0}^{t} f(t-\tau) g(\tau) d \tau \\
&=\mathcal{L}^l\left\{\int_{0}^{t} f(t-\tau) g(\tau) d \tau\right\} \\
&=\mathcal{L}^l\{f(t) \circ g(t)\}.
\end{aligned}
$$
Thus, $\mathcal{L}^l\{f\circ g\}$ is the left regular function on $U$, such that for $s\in \mathbb{R}$  and $s>c$, $\mathcal{L}^l\{f\circ g\}(s)=F(s) G(s).$ As $F(s)$ and $G(s)$ are in $\mathcal{R}^l(U)$, $\mathcal{L}^l\{f\circ g\}(s)=(F*G)(s)$ for $s\in U$ by Proposition \ref{k} (1).
\ep

\begin{prop} (Derivatives of the Laplace Transform).  Assume that $f(t)$ is of exponential order $a$ and  $U=\{s\in \mathbb{H}| \operatorname{Re}(s)>a\}$.  If $\mathcal{L}^l\{f(t)\}=F(s),$ then for $n\in \mathbb{N}$ and $s\in U$, 
$$
\begin{aligned}
\mathcal{L}^l\left\{t^{n} f(t)\right\} &=(-1)^{n}  F^{(n)}(s),
\end{aligned}
$$
where $F^{(n)}(s)$ is the $n$-th left slice derivative of $F(s)$.
\end{prop}
This can be proved as in the complex case.

\begin{prop}(The Laplace Transform of an Integral).  Assume that $f(t)$ is of exponential order $a$ and  $U=\{s\in \mathbb{H}| \operatorname{Re}(s)>a\}$. If $\mathcal{L}^l\{f(t)\}=F(s),$ then for $s\in U$,
$$
\mathcal{L}^l\left\{\int_{0}^{t} f(\tau) d \tau\right\}=s^{-1}{F(s)}.
$$
\end{prop}
This can also be proved as in the complex case using Proposition \ref{x}, see THEOREM 3.6.4 of \cite{it}.

So we have already seen that most of the properties for the complex Laplace transform can be generalized to the quaternionic Laplace transform.

\end{document}